\documentclass{elsart_mod}

\usepackage{natbib}

\usepackage[dvips]{graphicx}
\usepackage{subfigure}
\usepackage{psfrag}
\usepackage{amsmath}
\usepackage{amssymb}
\usepackage{enumerate}

\newcommand{\eref}[1]{(\ref{#1})}

\begin{document}
\begin{frontmatter}
\title{Characterizing the generalized lambda distribution by L-moments}
\author{Juha Karvanen\corauthref{cor}}
\address{Department of Health Promotion and Chronic Disease Prevention, National Public Health Institute,
Mannerheimintie 166, 00300 Helsinki, Finland}
\ead{juha.karvanen@ktl.fi}
\corauth[cor]{Corresponding author. Tel. +358 9 4744 8641, Fax +358 9 4744 8338}
\author{Arto Nuutinen}
\address{University of Helsinki, V{\"a}in{\"o} Auerin katu 1, 00560 Helsinki, Finland}
\maketitle

\begin{abstract}
The generalized lambda distribution (GLD) is a flexible four parameter distribution with many practical applications. L-moments of the GLD can be expressed in closed form and are good alternatives for the central moments. The L-moments of the GLD up to an arbitrary order are presented, and a study of L-skewness and L-kurtosis that can be achieved by the GLD is provided. The boundaries of L-skewness and L-kurtosis are derived analytically for the symmetric GLD and calculated numerically for the GLD in general. Additionally, the contours of L-skewness and L-kurtosis are presented as functions of the GLD parameters. It is found that with an exception of the smallest values of L-kurtosis, the GLD covers all possible pairs of L-skewness and L-kurtosis and often there are two or more distributions that share the same L-skewness and the same L-kurtosis. Examples that demonstrate situations where there are four GLD members with the same L-skewness and the same L-kurtosis are presented. The estimation of the GLD parameters is studied in a simulation example where method of L-moments compares favorably to more complicated estimation methods. The results increase the knowledge on the distributions that belong to the GLD family and can be utilized in model selection and estimation.
\end{abstract}
\end{frontmatter}

\section{Introduction}
The generalized lambda distribution (GLD) is a four parameter distribution that has been applied to various problems where a flexible parametric model for univariate data is needed. The GLD provides fit for a large range of skewness and kurtosis values and can approximate many commonly used distributions such as normal, exponential and uniform. The applications of the GLD include e.g. option pricing \citep{Corrado:optionpricing}, independent component analysis \citep{vlsi}, statistical process control \citep{Pal:processgld}, analysis of fatigue of materials \citep{Bigerelle:fatigue}, measurement technology \citep{Lampasi:glduncertainty} and generation of random variables \citep{Ramberg:74,gldrnd,Headrick:simulating}.

\begin{figure}[htb]
\begin{center}
\psfrag{lm4}{$\lambda_4$}
  \psfrag{lm3}{$\lambda_3$}
\includegraphics[width=0.6\textwidth]{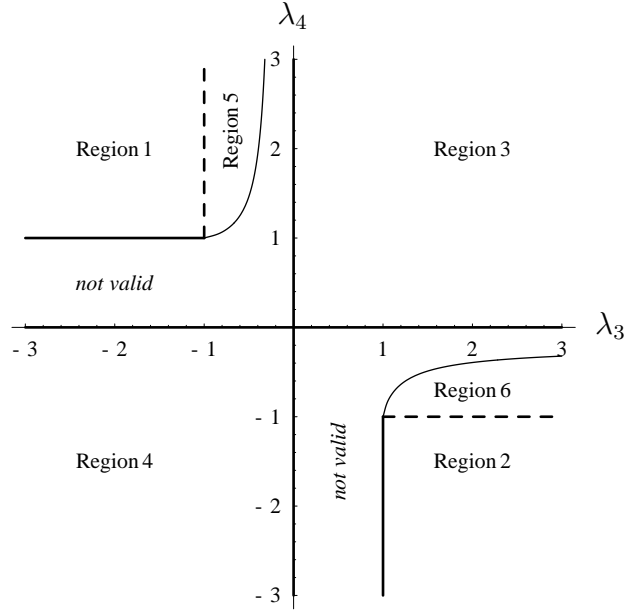}
\end{center}
\label{fig:gldregions}
\caption{GLD regions as defined in \citep{egld:theory}. The distributions are bounded on the right in regions 1 and 5, bounded on the left in regions 2 and 6, bounded in region 3 and unbounded in region 4.}
\end{figure}

The price for the high flexibility is high complexity. Probability density function (pdf) or cumulative distribution function (cdf) of the GLD do not exist in closed form but the distribution is defined by the inverse distribution function \citep{Ramberg:74}
\begin{equation}
\label{eq:glddefinition}
F^{-1}(u)=\lambda_1+\frac{u^{\lambda_3}-(1-u)^{\lambda_4}}{\lambda_2},
\end{equation}
where $0 \leq u \leq 1$ and $\lambda_1$, $\lambda_2$, $\lambda_3$ and $\lambda_4$ are the parameters of the GLD.  Equation~\eref{eq:glddefinition} defines a distribution if and only if \citep{egld:theory}
\begin{equation}
\frac{\lambda_{2}}{\lambda_{3}u^{\lambda_{3}-1}+\lambda_{4}(1-u)^{\lambda_{4}-1}} \geq 0 \textrm{ for all } u \in [0,1].
\end{equation}
\citet{egld:theory} divide the GLD into six regions on the basis of the parameters $\lambda_3$ and $\lambda_4$ that control skewness and kurtosis of the distribution. These regions are presented in Figure~\ref{fig:gldregions}. The regions have different characteristics: the distributions in region 3 are bounded whereas the distributions in region 4 are unbounded; the distributions in regions 1 and 5 are bounded on the right and the distributions in regions 2 and 6 are bounded on the left. The boundaries of the domain in each region are described e.g. by \citet{KarianDudewicz:gldbook} and \citet{Fournier:estimating}.

The members of the GLD family have been traditionally characterized by the central moments, especially by the central moment skewness and kurtosis \citep{Ramberg:74,egld:theory}. L-moments \citep{Hosking:Lmoments}, defined as linear combinations of order statistics, are attractive alternatives for the central moments. Differently from the central moments, the L-moments of the GLD can be expressed in closed form, which allows us to derive some analytical results and makes it straightforward to perform numerical analysis. The parameters of the GLD can be estimated by method of L-moments \citep{vlsi,Asquith:lmomgld}. Method of L-moments can be used independently or together with other estimation methods, analogously to what was done in \citep{Su:numericalmaximum}, \citep{Fournier:estimating} and \citep{Lakhany:estimating} with numerical maximum likelihood, method of moments, method of percentiles \citep{egld:percentile}, the least square method \citep{Ozturk:leastsquares} and the starship method \citep{King:starsip}.

\citet{Karian:comparison} report that method of percentiles gives superior fits compared to method of moments. They also study method of L-moments in some of their examples where the overall performance looks comparable to the overall performance of the percentile method. Method of percentiles and method of L-moments are related in sense that they both are based on order statistics. There are, however, some differences: L-moments are linear combinations of all order statistics whereas method of percentiles uses only a limited number of order statistics that need to be explicitly chosen. \citet{Fournier:estimating} study the choice of the order statistics for the percentile method and conclude that there is no trivial rule for choosing them. Such problems are avoided in method of L-moments.

In this paper we analyze the GLD using L-moments and consider the estimation by method of L-moments. In Section~\ref{sec:gldlmoments} we present the L-moments of the GLD up to an arbitrary order. In Section~\ref{sec:symmetric} we analyze the symmetric case $(\lambda_{3}=\lambda_{4})$ and derive the boundaries of L-kurtosis analytically and in Section~\ref{sec:boundaries} we consider the general case $(\lambda_{3} \neq \lambda_{4})$ and calculate the boundaries of the GLD in the terms of L-skewness and L-kurtosis using numerical methods. The boundaries are calculated separately for each GLD region. We also calculate the contours of L-skewness and L-kurtosis as functions of $\lambda_3$ and $\lambda_4$. The examples presented illustrate the multiplicity of the distributional forms of the GLD. In Section~\ref{sec:estimation} we present an example on the estimation by method of L-moments and compare the results to the results by \citet{Fournier:estimating}. Section~\ref{sec:conclusion} concludes the paper.

\section{L-moments and the GLD} \label{sec:gldlmoments}
The L-moment of order $r$ can be expressed as \citep{Hosking:Lmoments}
\begin{equation}
L_{r}=\int_{0}^{1} F^{-1}(u) P_{r-1}^{*}(u) \textrm{d}u,
\end{equation}
where
\begin{equation}
P_{r-1}^{*}(u)=\sum_{k=0}^{r-1} (-1)^{r-k-1} \binom{r-1}{k} \binom{r+k-1}{k} u^{k}
\end{equation}
is the shifted Legendre polynomial of order $r-1$. All L-moments of a real-valued random variable exists if and only if the random variable has a finite mean and furthermore, a distribution whose mean exists, is uniquely determined by its L-moments \citep{Hosking:Lmoments,Hosking:characterization}.

Similarly to the central moments, $L_1$ measures location and $L_2$ measures scale. The higher order L-moments are usually transformed to L-moment ratios
\begin{equation}
\tau_r = \frac{L_{r}}{L_2} \quad r=3,4,\ldots
\end{equation}
L-skewness $\tau_{3}$ is related to the asymmetry of the distribution and L-kurtosis $\tau_{4}$ is related to the peakedness of the distribution. Differently from the central moment skewness and kurtosis, $\tau_{3}$ and $\tau_{4}$ are constrained by the conditions \citep{Hosking:Lmoments, Jones:someexpressions}
\begin{align}
-1 < &\tau_{3} < 1 \\
\intertext{and}
(5\tau_{3}^{2}-1)/4 \leq &\tau_{4} < 1.
\end{align}
Distributions are commonly characterized using an L-moment ratio diagram in which L-skewness is on the horizontal axis and L-kurtosis is on the vertical axis.

The L-moments of the GLD have been presented in the literature up to order~5 \citep{Bergevin:analysis,vlsi,Asquith:lmomgld}. We generalize the results for an arbitrary order $r$. The derivation is straightforward if we first note that
\begin{align}
\int_{0}^{1} u^{\lambda_{3}} P_{r-1}^{*}(u) \textrm{d}u=&\sum_{k=0}^{r-1} \frac{(-1)^{r-k-1} \binom{r-1}{k} \binom{r+k-1}{k}}{k+1+\lambda_{3}}\\
\intertext{and} \int_{0}^{1} (1-u)^{\lambda_{4}} P_{r-1}^{*}(u) \textrm{d}u=&-\int_{0}^{1} (u)^{\lambda_{4}} P_{r-1}^{*}(1-u)\textrm{d}u \nonumber \\ =&(-1)^{r}\sum_{k=0}^{r-1} \frac{(-1)^{r-k-1} \binom{r-1}{k} \binom{r+k-1}{k}}{k+1+\lambda_{4}}, \label{eq:ld4u}
\end{align}
where the last equality follows from the property $P_{r-1}^{*}(1-u)=(-1)^{r-1}P_{r-1}^{*}(u)$.
We obtain for $r \geq 2$
\begin{equation}
\lambda_{2}L_{r}=\sum_{k=0}^{r-1} (-1)^{r-k-1} \binom{r-1}{k} \binom{r+k-1}{k} \left(\frac{1}{k+1+\lambda_{3}}+\frac{(-1)^{r}}{k+1+\lambda_{4}} \right).
\end{equation}
The explicit formulas for the first six L-moments of the GLD are
\begin{align}
L_{1}=&\lambda_1-\frac{1}{\lambda_2} \left( \frac{1}{1+\lambda_4} - \frac{1}{1+\lambda_3} \right), \label{eq:L1}\\
L_2 \lambda_2 =& - \frac{1}{1 + \lambda_3} + \frac{2}{2 + \lambda_3} - \frac{1}{1 + \lambda_4} + \frac{2}{2 + \lambda_4}, \label{eq:L2} \\
L_3 \lambda_2 =& \frac{1}{1+\lambda_3} - \frac{6}{2+\lambda_3} + \frac{6}{3+\lambda_3} - \frac{1}{1+\lambda_4} + \frac{6}{2+\lambda_4} - \frac{6}{3+\lambda_4}, \label{eq:L3}\\
L_4 \lambda_2 =& - \frac{1}{1+\lambda_3} + \frac{12}{2+\lambda_3} - \frac{30}{3+\lambda_3} + \frac{20}{4+\lambda_3}\nonumber \\
&-\frac{1}{1+\lambda_4} +  \frac{12}{2+\lambda_4} - \frac{30}{3+\lambda_4} + \frac{20}{4+\lambda_4}, \label{eq:L4}\\
L_5 \lambda_2 =&\frac{1}{1+\lambda_3} - \frac{20}{2+\lambda_3} + \frac{90}{3+\lambda_3} -\frac{140}{4+\lambda_3}+\frac{70}{5+\lambda_3}\nonumber \\
-&\frac{1}{1+\lambda_4} + \frac{20}{2+\lambda_4} - \frac{90}{3+\lambda_4} +\frac{140}{4+\lambda_4}-\frac{70}{5+\lambda_4},\\
\intertext{and} L_6 \lambda_2 =& -\frac{1}{1+\lambda_3} + \frac{30}{2+\lambda_3} - \frac{210}{3+\lambda_3} +\frac{560}{4+\lambda_3}-\frac{630}{5+\lambda_3}+\frac{252}{6+\lambda_3}\nonumber \\
&-\frac{1}{1+\lambda_4} + \frac{30}{2+\lambda_4} - \frac{210}{3+\lambda_4} +\frac{560}{4+\lambda_4}-\frac{630}{5+\lambda_4}+\frac{252}{6+\lambda_4}.
\end{align}
Note that the results in \citep{Asquith:lmomgld} and \citep{Bergevin:analysis} are the same but expressed in a different form.

The mean of GLD and therefore also all L-moments exist if $\lambda_{3},\lambda_{4}>-1$. This implies that the characterization by L-moments covers regions 3, 5 and 6 and the subset of region 4 where $\lambda_{3},\lambda_{4}>-1$.

\section{The L-moments of the GLD in the symmetric case} \label{sec:symmetric}
\begin{figure}[htb]
\begin{center}
\psfrag{t}{$\! \! \! \tau_4$}
  \psfrag{l}{$\lambda_3$}
\includegraphics[width=0.8\textwidth]{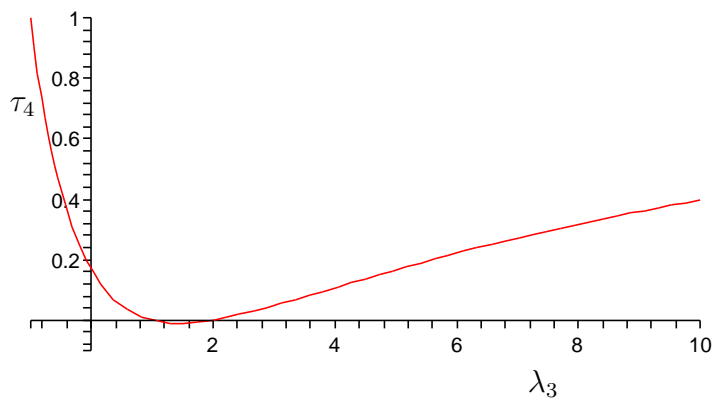}
\end{center}
\caption{ L-kurtosis $\tau_4(\lambda_{3})$ for symmetric distributions $\lambda_3=\lambda_4$. \label{fig:symmtau4}}
\end{figure}

In this section we consider the special case $\lambda_{3}=\lambda_{4}$, which defines a symmetric distribution. The symmetric distributions can be bounded (in region 3) or unbounded (in region 4). Applying the condition $\lambda_{3}=\lambda_{4}$ to equations \eref{eq:L3} and \eref{eq:L4} we find that $\tau_{3}=0$ and
\begin{equation} \label{eq:symmtau4}
\tau_{4}=\frac{\lambda_{3}^2-3\lambda_{3}+2}{\lambda_{3}^2+7\lambda_{3}+12}, \quad \lambda_{3}=\lambda_{4} > -1.
\end{equation}
In Figure \ref{fig:symmtau4} $\tau_{4}$ is plotted as a function of $\lambda_{3}$. The value of $\tau_{4}$ approaches to 1 as $\lambda_{3}$ approaches to infinity. The minimum of $\tau_{4}$ is found to be
\begin{equation}
\frac{12-5\sqrt{6}}{12+5\sqrt{6}} \approx -0.0102,
\end{equation}
and is obtained when \mbox{$\lambda_{3}=-1+\sqrt{6}$}.

In the symmetric case it is also possible to solve $\lambda_{3}$ and $\lambda_{4}$ as functions of $\tau_{4}$ \citep{vlsi}
\begin{equation} \label{eq:symmsolution}
\lambda_4 =\lambda_3=\frac{3+7\tau_4 \pm \sqrt{1+98\tau_4+\tau_4^2}}{2(1-\tau_4)}.
\end{equation}
It follows that if $\tau_{4}>1/6$ there are symmetric GLD members available from both regions 3 and 4. If
\begin{equation}
\frac{12-5\sqrt{6}}{12+5\sqrt{6}} < \tau_{4} \leq \frac{1}{6},
\end{equation}
there are two GLD members available from region 3, and if
\begin{equation}
\tau_{4} < \frac{12-5\sqrt{6}}{12+5\sqrt{6}},
\end{equation}
there are no symmetric GLD members available.

\begin{figure}[htb]
\begin{center}
\subfigure[$\tau_{4} \approx 0.12260$ ($\tau_{4}$ of normal distribution)]{
 \includegraphics[width=0.48\columnwidth]{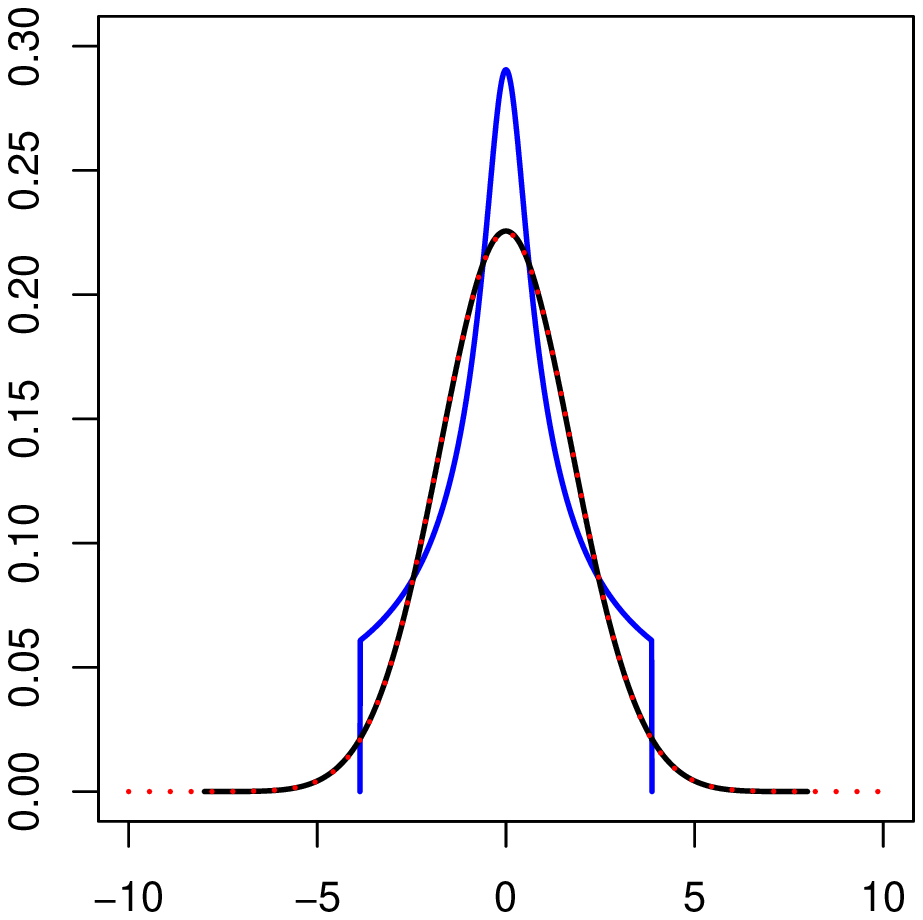}}
\subfigure[$\tau_{4}=0.25$]{
 \includegraphics[width=0.48\columnwidth]{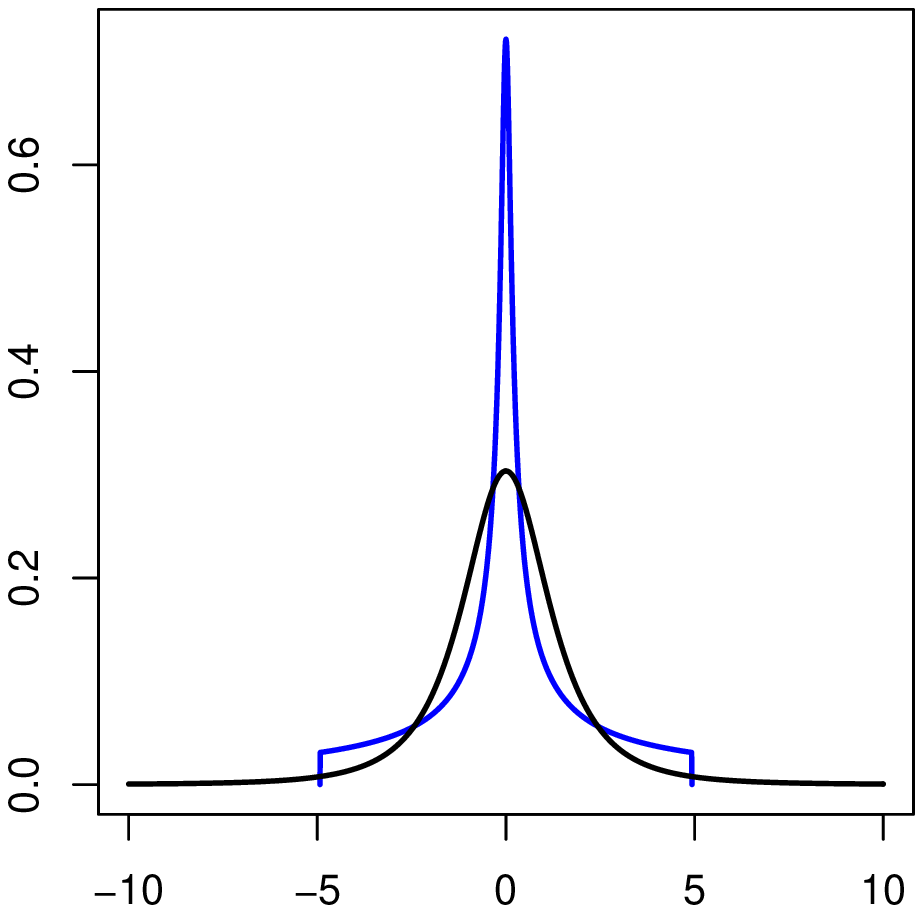}}
\end{center}
\caption{Examples of symmetric GLDs sharing the same $\tau_{4}$. The pdfs of the GLDs are plotted with solid line and the pdf of normal distribution is plotted with dotted line. All distributions have $\lambda_{1}=0$, $\lambda_{2}=1$ and $\tau_{3}=0$.} \label{fig:symmpdfs}
\end{figure}

Examples of symmetric GLDs sharing the same $\tau_{4}$ are presented in Figure~\ref{fig:symmpdfs}. In Figure~\ref{fig:symmpdfs}(a) the value of $\tau_{4}$ is set to be equal to the $\tau_{4}$ of the normal distribution. Both GLDs are from region 3 and are bounded.  The GLD for the more peaked distribution in Figure~\ref{fig:symmpdfs}(a) has $\tau_{6} \approx 0.0004$ whereas the other GLD has $\tau_{6} \approx 0.043$, which is very near to the $\tau_{6}$ of the normal distribution. For comparison the pdf of the normal distribution is plotted (dotted line). The differences with the normal distribution are appreciable only in the far tails. In Figure~\ref{fig:symmpdfs}(b) $\tau_{4}=0.25$, which means that equation~\eref{eq:symmtau4} has one positive and one negative root. The GLD with $\tau_{6} \approx 0.016$ and  the higher peak is from region 3 and has a bounded domain. The other GLD with $\tau_{6} \approx 0.121$ is from region 4 and has unbounded domain.

In addition to the symmetric case we may consider some other special cases. If $\lambda_{3}=0$ we obtain
\begin{align}
\tau_{3}(\lambda_{4})=&\frac{1-\lambda_{4}}{\lambda_{4}+3} \\
\intertext{and} \tau_{4}(\lambda_{4})=&\frac{\lambda_{4}^{2}-3\lambda_{4}+2}{\lambda_{4}^{2}+7\lambda_{4}+12}. \label{eq:tau4ld4}
\end{align}
If  $\lambda_{4}=0$ we obtain
\begin{align}
\tau_{3}(\lambda_{3})=&\frac{\lambda_{3}-1}{\lambda_{3}+3} \\
\intertext{and} \tau_{4}(\lambda_{3})=&\frac{\lambda_{3}^{2}-3\lambda_{3}+2}{\lambda_{3}^{2}+7\lambda_{3}+12}. \label{eq:tau4ld3}
\end{align}
Interestingly, equations~\eref{eq:tau4ld4} and \eref{eq:tau4ld3} are identical to equation~\eref{eq:symmtau4}, which indicates that the $(\lambda_{3},\lambda_{4})$ pairs $(\lambda,\lambda)$, $(\lambda,0)$ and $(0,\lambda)$ lead to the same value of $\tau_{4}$.

\section{Boundaries of the L-moment ratios of the GLD} \label{sec:boundaries}
In the general case it is difficult to derive analytical results and therefore we resort to numerical methods. Because the L-moments of the GLD are available in closed form, we may choose a straightforward approach and calculate the values of $\tau_{3}$ and $\tau_{4}$ for a large number of $(\lambda_{3},\lambda_{4})$ pairs. Naturally, the values of $(\lambda_{3},\lambda_{4})$ need to be chosen such a way that the essential properties of the GLD are revealed. In region 4, we use a grid of one million points where both $\lambda_{3}$ and $\lambda_{4}$ are equally spaced in the interval $[-1+10^{-10},10^{-10}]$. In region 3 we use an unequally spaced grid of 4104676 points. The grid is dense near the origin ($\lambda_{3}$ and $\lambda_{4}$ have minimum $10^{-10}$) and sparse for large values ($\lambda_{3}$ and $\lambda_{4}$ have maximum $10^{12}$). In region 5, $\lambda_{3}$ is equally spaced in the interval $[-1+10^{-10},0]$ and $\lambda_{4}$ is unequally spaced starting from $-1$. The grid contains 1362429 valid points. The calculations thus result in a large number of $(\tau_{3},\tau_{4})$ pairs. The results for region 6 can be obtained from the results for region 5 by swapping $\lambda_{3}$ and $\lambda_{4}$. The problem is to find the boundaries of the GLD in the $(\tau_{3},\tau_{4})$ space on the basis of these data. Let $(\tau_{3}[i],\tau_{4}[j])$ be the image of the grid point $(\lambda_{3}[i],\lambda_{4}[j])$, where $i$ and $j$ are the grid indexes. If the point $(\tau_{3}[i],\tau_{4}[j])$ is not located inside the polygon defined by points $(\tau_{3}[i-1],\tau_{4}[j])$, $(\tau_{3}[i],\tau_{4}[j+1])$, $(\tau_{3}[i+1],\tau_{4}[j])$ and $(\tau_{3}[i],\tau_{4}[j-1])$ it is considered to be a potential boundary point. The images of the points that are located on the boundary of the grid in the $(\lambda_{3},\lambda_{4})$ space are also potential boundary points in the $(\tau_{3},\tau_{4})$ space. The actual boundaries are found combining the potential boundaries and finally it is checked that all $(\tau_{3},\tau_{4})$ pairs are located inside the boundaries. The calculations are made using R \citep{R} and the contributed R packages \texttt{PBSmapping} \citep{R_PBSmapping} and \texttt{gld} \citep{R_gld} are utilized.

\begin{figure}[htb]
\begin{center}
\subfigure[Region 3]{
 \includegraphics[width=0.48\columnwidth]{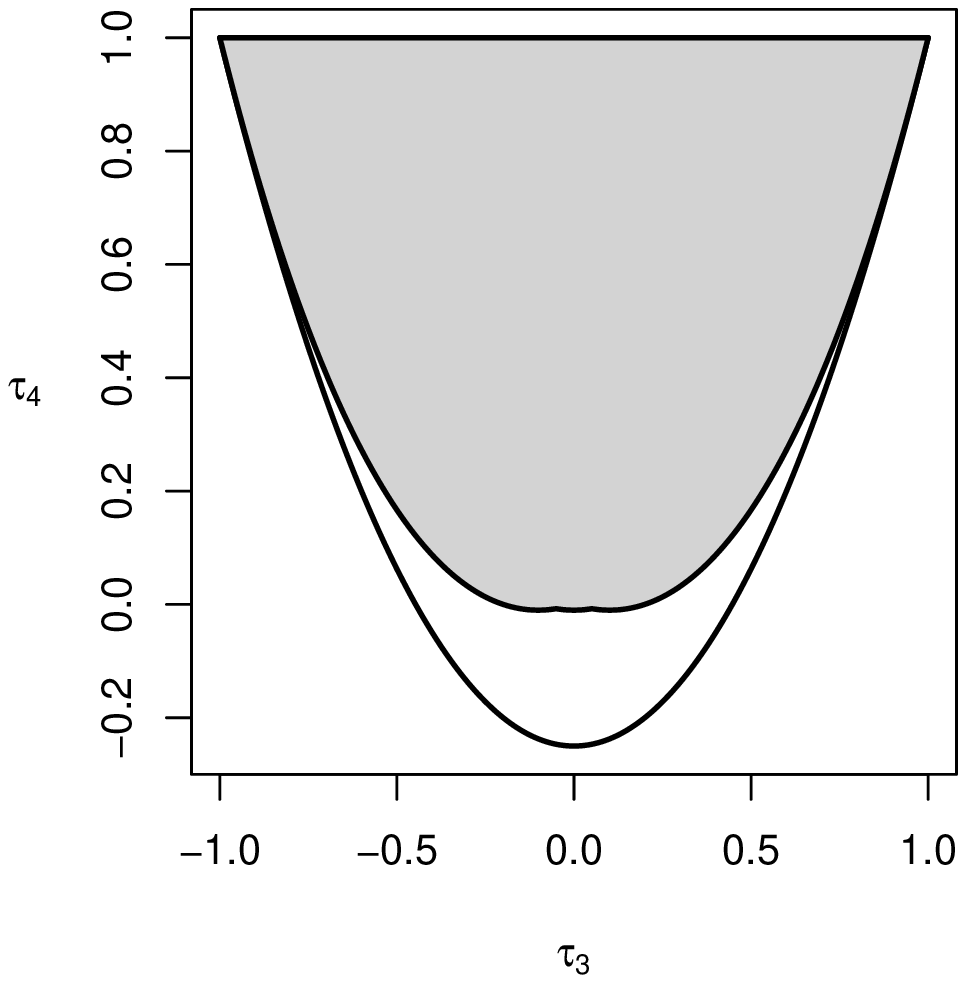}}
\subfigure[Region 4]{
 \includegraphics[width=0.48\columnwidth]{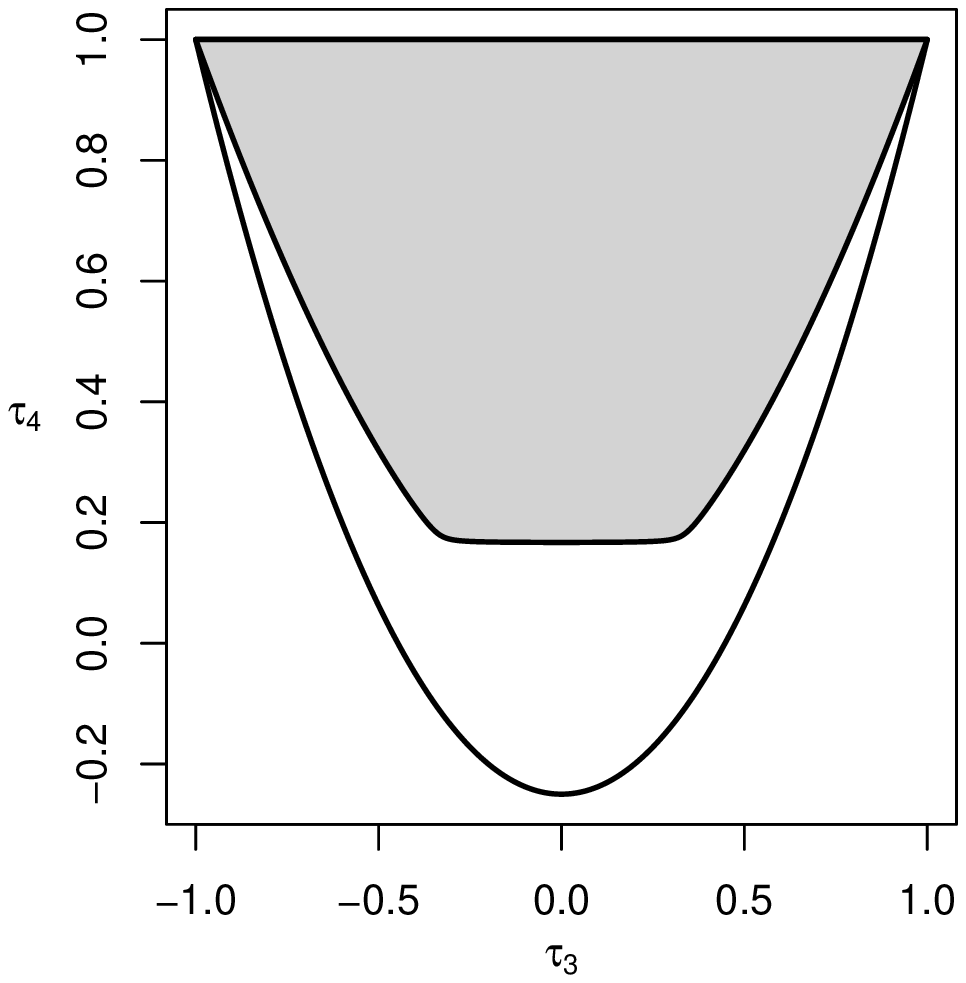}}
 \subfigure[Region 5]{
 \includegraphics[width=0.48\columnwidth]{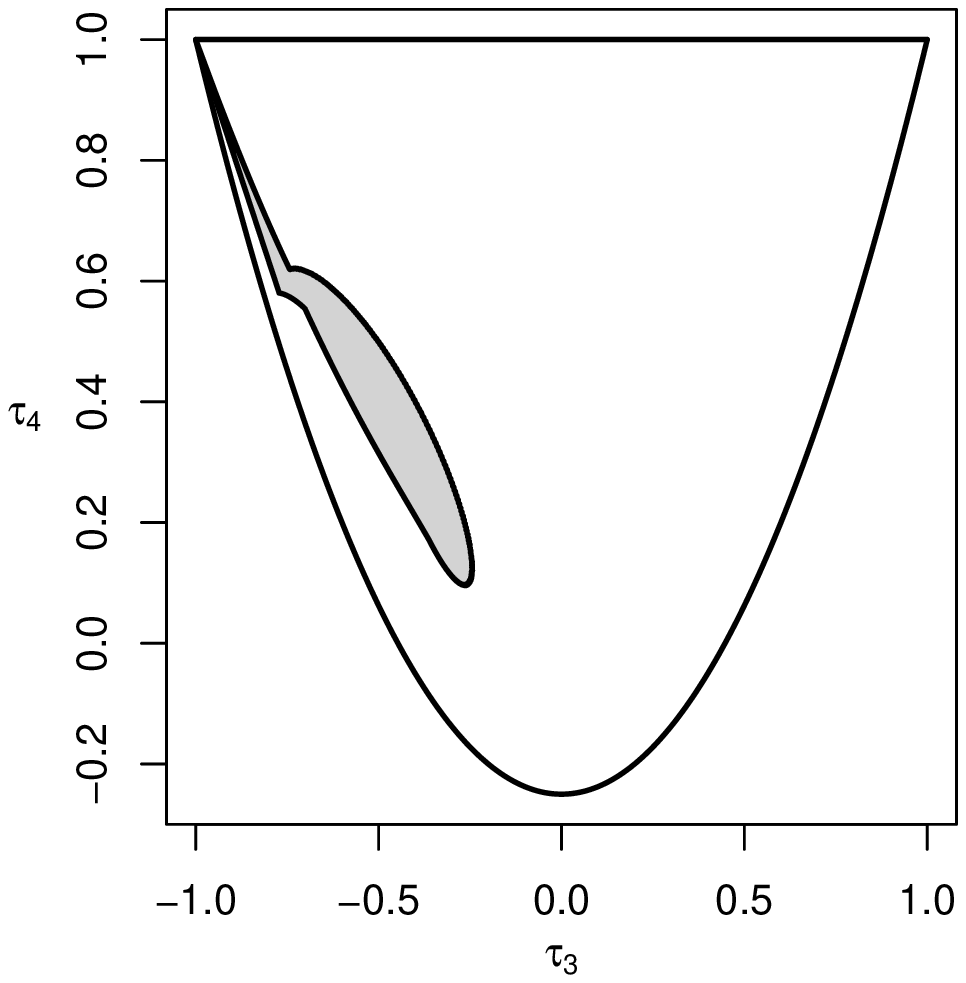}}
\subfigure[Region 6]{
 \includegraphics[width=0.48\columnwidth]{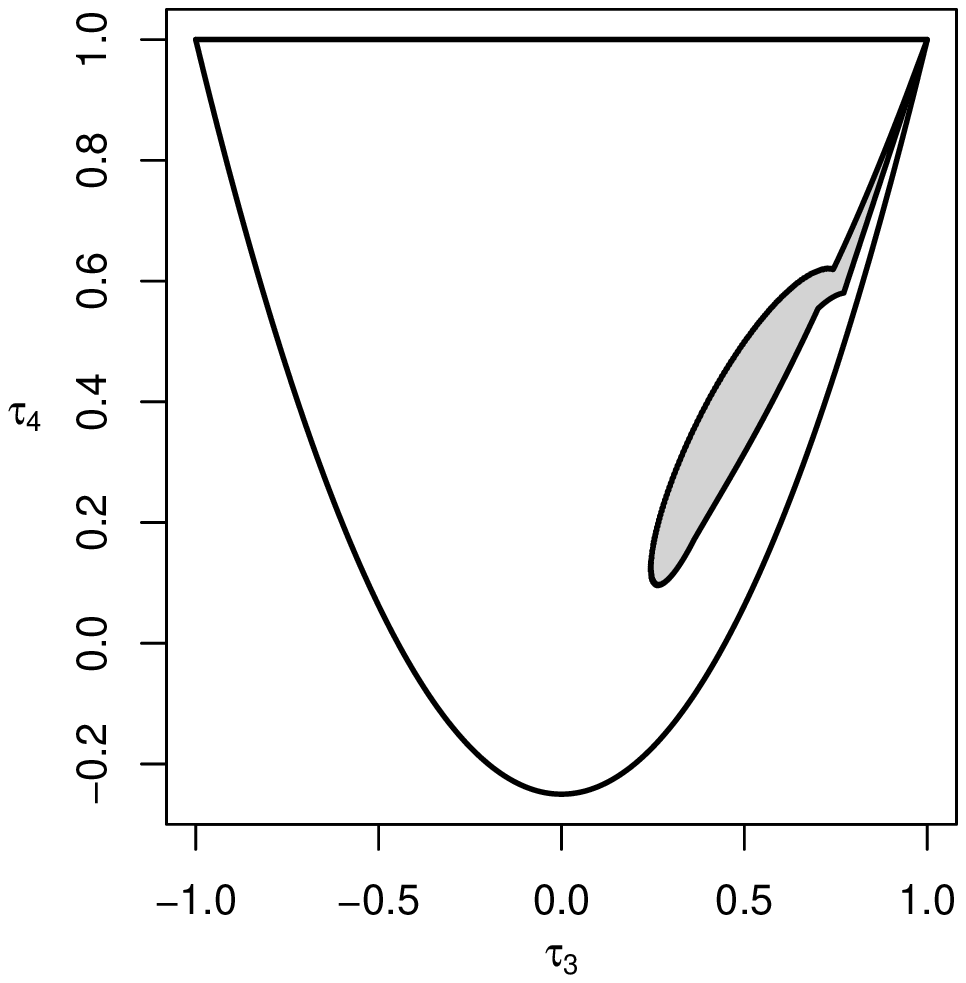}}
\end{center}
\caption{Boundaries of the GLD in the $(\tau_{3},\tau_{4})$ space. The shading shows the values of L-skewness and L-kurtosis that are achievable by the GLD. The outer limits present the boundaries for all distributions.} \label{fig:boundaries}
\end{figure}

The boundaries are shown in Figure~\ref{fig:boundaries}. It can be seen that region 3 covers most of the $(\tau_{3},\tau_{4})$ space; only the smallest values of $\tau_{4}$ are unattainable. Region 4 largely overlaps with region 3 but cannot achieve $\tau_{4}$ smaller than $1/6$. Region 5 and its counterpart, region 6, cover a rather small area of the $(\tau_{3},\tau_{4})$ space.

\begin{figure}[htb]
\begin{center}
\includegraphics[width=0.48\columnwidth]{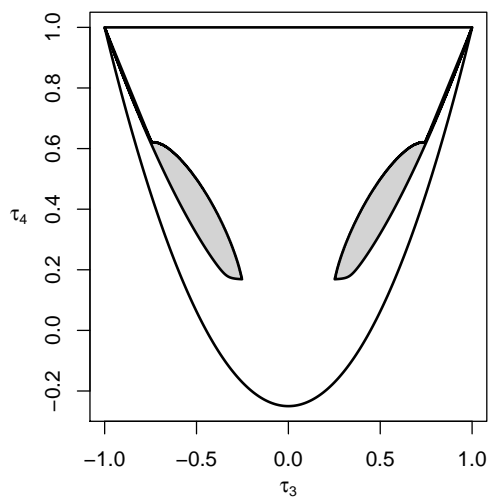}
\end{center}
\caption{The $(\tau_{3},\tau_{4})$ area where there exist GLD members from regions 3, 4 and 5 (negative L-skewness), or from regions 3, 4 and 6 (positive L-skewness).} \label{fig:reg3456int}
\end{figure}

\begin{figure}[htb]
\begin{center}
\subfigure[from region 3]{
 \includegraphics[width=0.48\columnwidth]{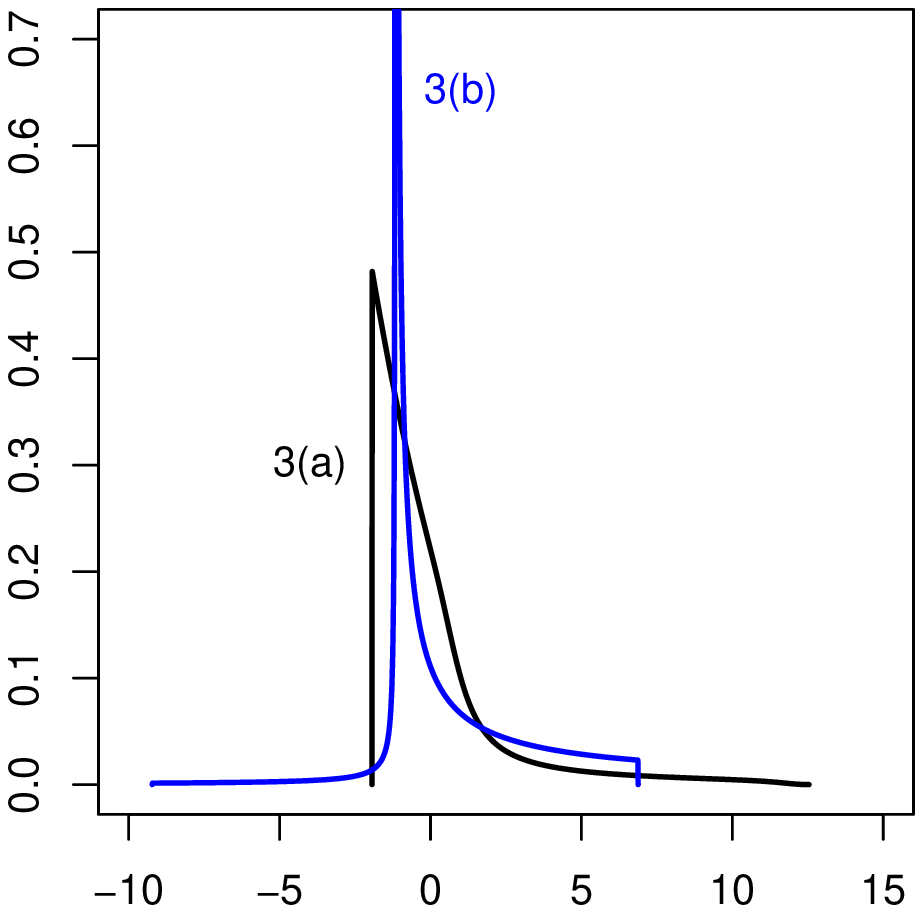}}
\subfigure[from regions 4 and 6]{
 \includegraphics[width=0.48\columnwidth]{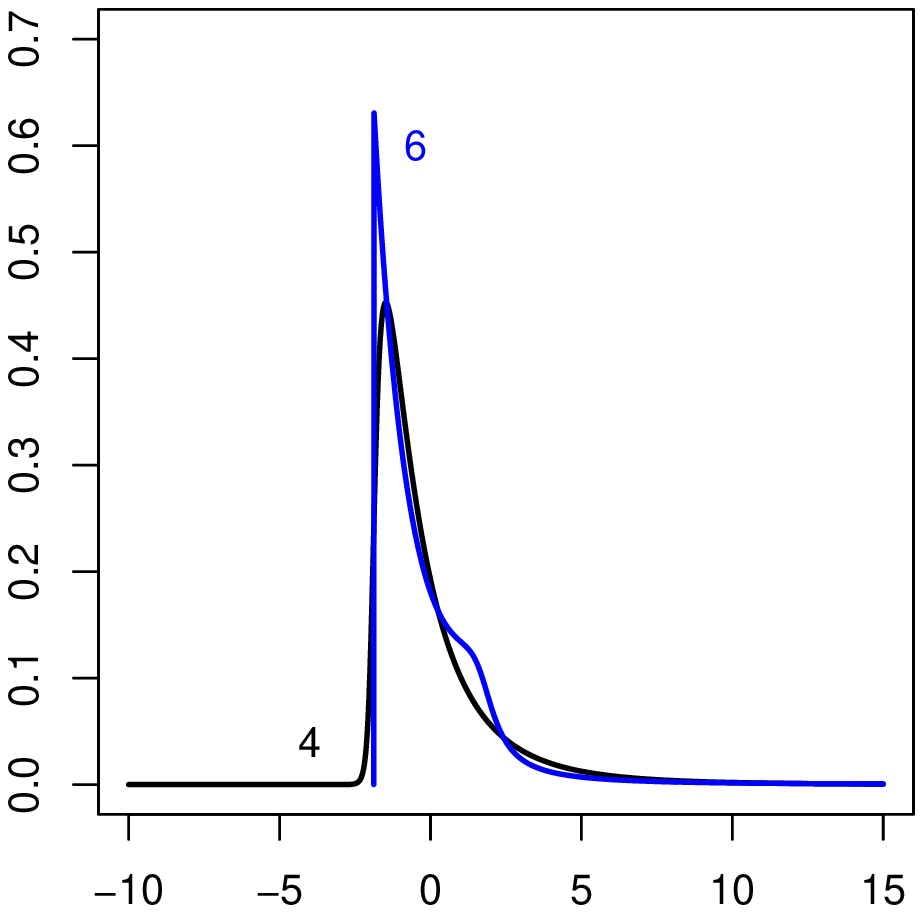}}
\end{center}
\caption{Pdfs of four GLDs with $\tau_{3}=0.4$ and $\tau_{4}=0.25$. The parameters of the distributions are presented in Table~\ref{tab:4gldexample}.} \label{fig:4gldexample}
\end{figure}

Figure~\ref{fig:reg3456int} shows the $(\tau_{3},\tau_{4})$ area where there exist GLD members from regions 3, 4 and 5, or from regions 3, 4 and 6.  To illustrate these distributions we fix $\lambda_{1}=0$, $\lambda_{2}=1$, $\tau_{3}=0.4$ and $\tau_{4}=0.25$ and seek for the solutions in regions 3, 4 and 6. Table~\ref{tab:4gldexample} and Figure~\ref{fig:4gldexample} present the four distributions that are found. Although the first four L-moments of the distributions are the same, the differences of the pdfs are clearly visible. Distribution 3(a) has bounded domain with sharp left limit. Distribution 3(b) is also bounded but characterized by the high peak. Distribution 4 has unbounded domain but is otherwise almost similar to distribution 3(a). Distribution 6 is bounded from the left and characterized by a minor `peak' around $x=2$.

\begin{table}
\caption{Four different GLD members with $L_{1}=0$, $L_{2}=1$, $\tau_{3}=0.4$ and $\tau_{4}=0.25$. \label{tab:4gldexample}}
\begin{center}
\begin{tabular}{crrrrrr}
Region & $\lambda_{1}$ & $\lambda_{2}$ & $\lambda_{3}$ & $\lambda_{4}$ & $\tau_{5}$ & $\tau_{6}$\\
3(a) & 5.322  &  0.138  &  21.526  &  0.286  &  0.163  &  0.103 \\
3(b) & -1.168  &  0.124  &  5.417  &  92.608  &  -0.029  &  0.067 \\
4 &  -1.62  &  -0.157  &  -0.014  &  -0.212  &  0.158  &  0.121 \\
6 & -7.04  &  -0.194  &  11.905  &  -0.306  &  0.204  &  0.180
\end{tabular}
\end{center}
\end{table}

\begin{figure}[htb]
\begin{center}
 \subfigure[Contours of L-skewness $\tau_{3}$]{
 \includegraphics[width=0.48\columnwidth]{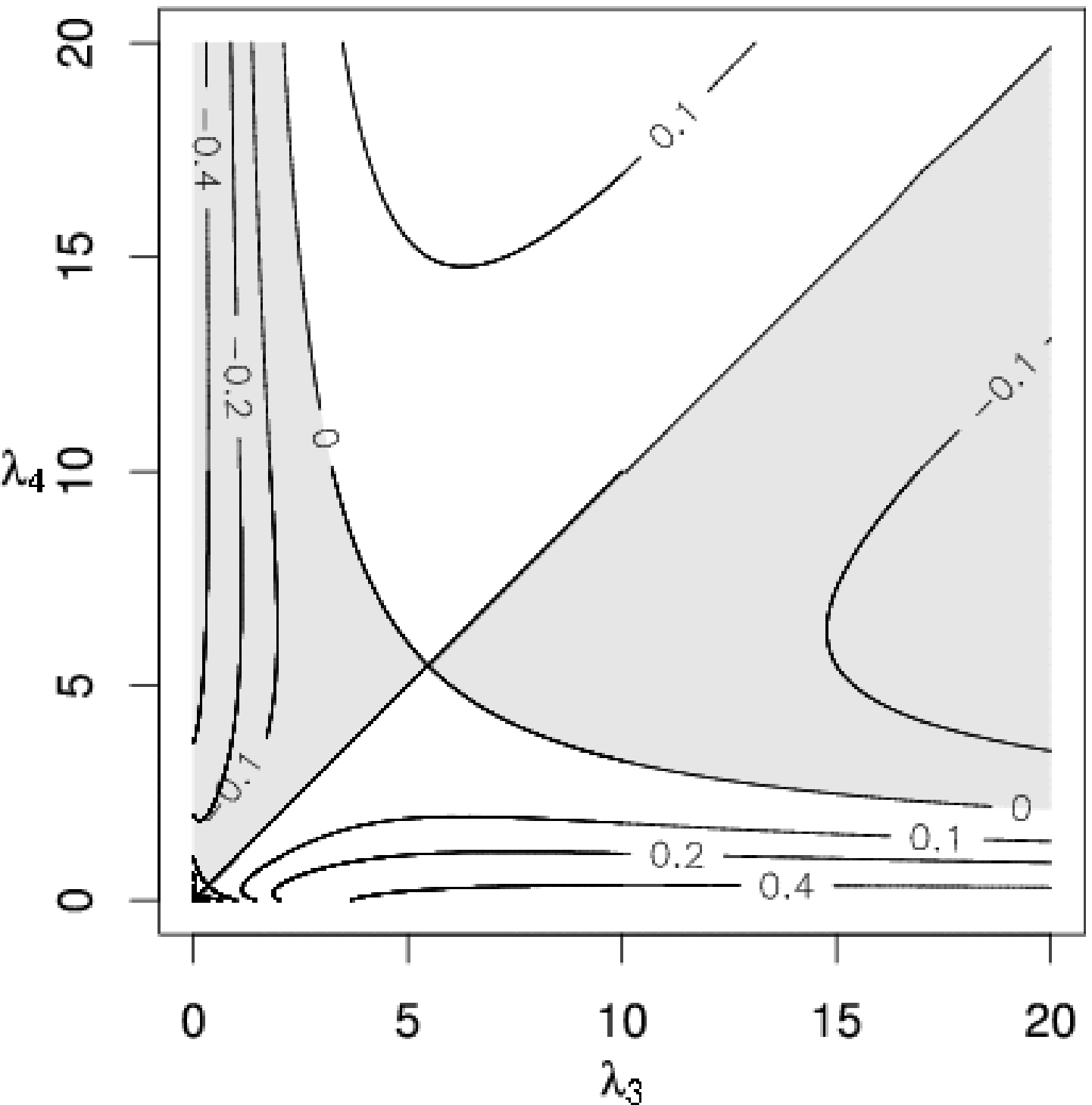}}
\subfigure[Contours of L-kurtosis $\tau_{4}$]{
 \includegraphics[width=0.48\columnwidth]{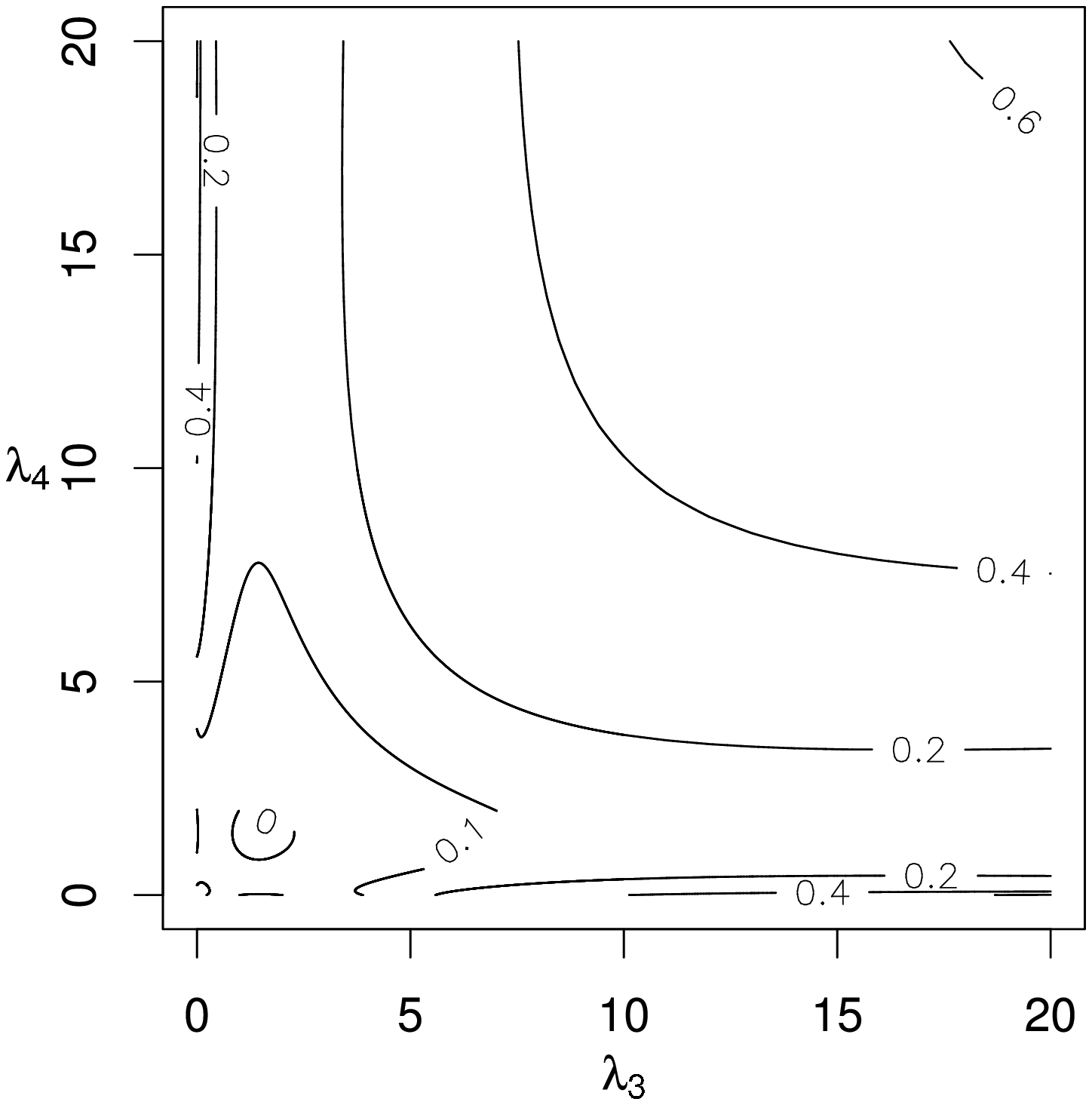}}
 \subfigure[Contours of L-skewness $\tau_{3}$ (zoomed)]{
 \includegraphics[width=0.48\columnwidth]{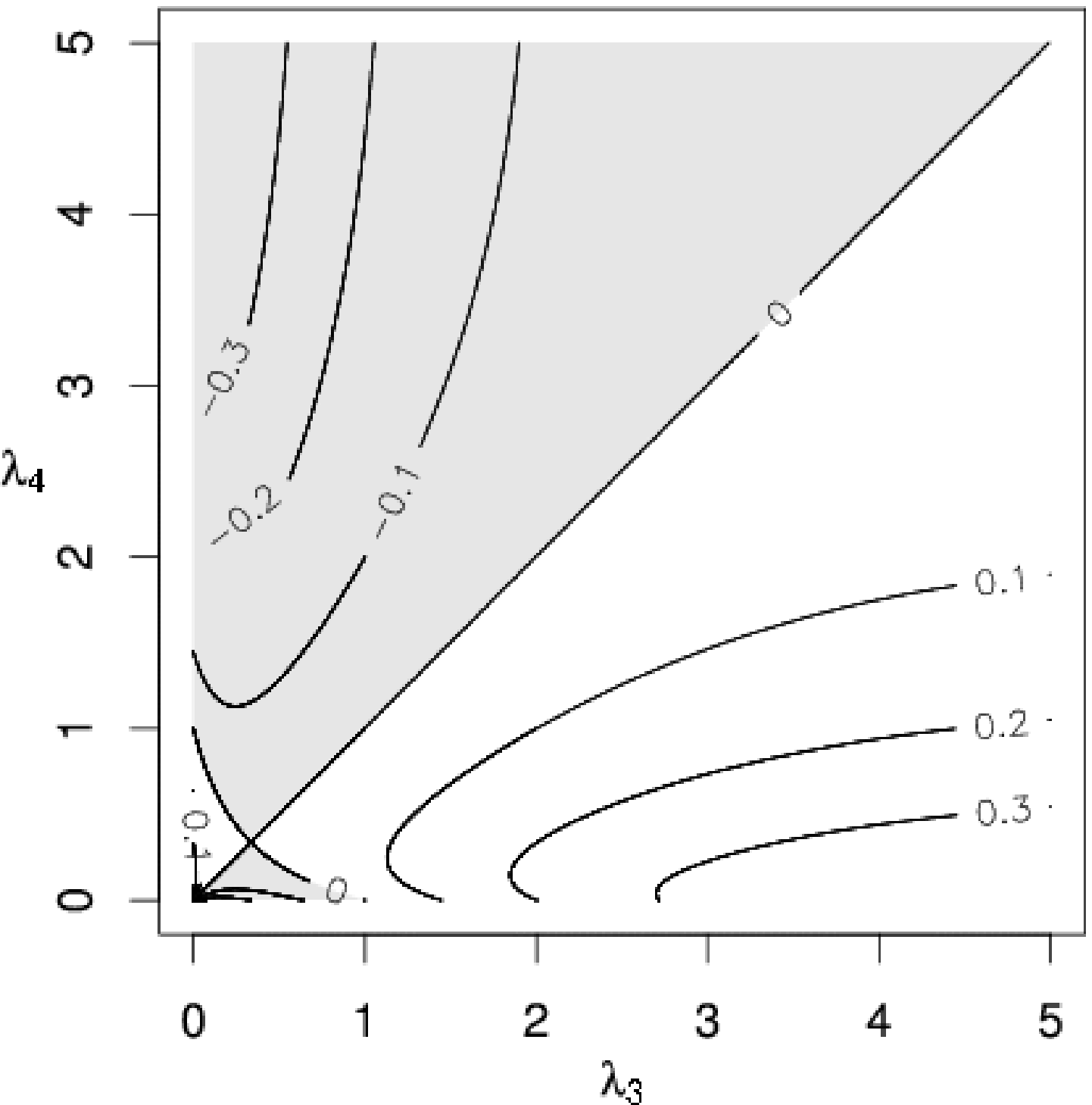}}
\subfigure[Contours of L-kurtosis $\tau_{4}$ (zoomed)]{
 \includegraphics[width=0.48\columnwidth]{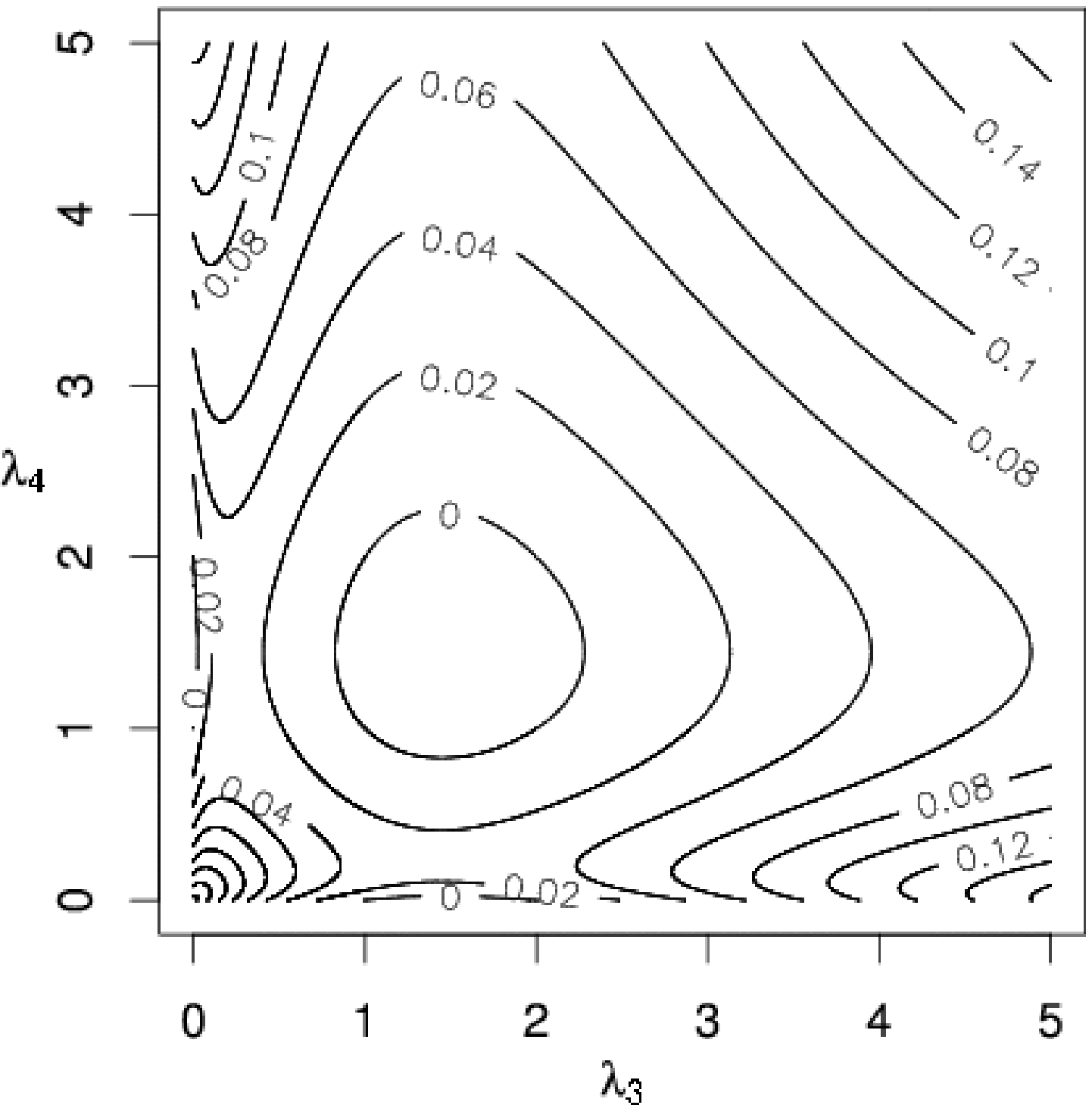}}
\end{center}
\caption{Contours of  $\tau_{3}$ and $\tau_{4}$ in region 3. The shading indicates negative values of $\tau_{3}$.} \label{fig:reg3contours}
\end{figure}

\begin{figure}[htb]
\begin{center}
\subfigure[Contours of L-skewness $\tau_{3}$]{
 \includegraphics[width=0.48\columnwidth]{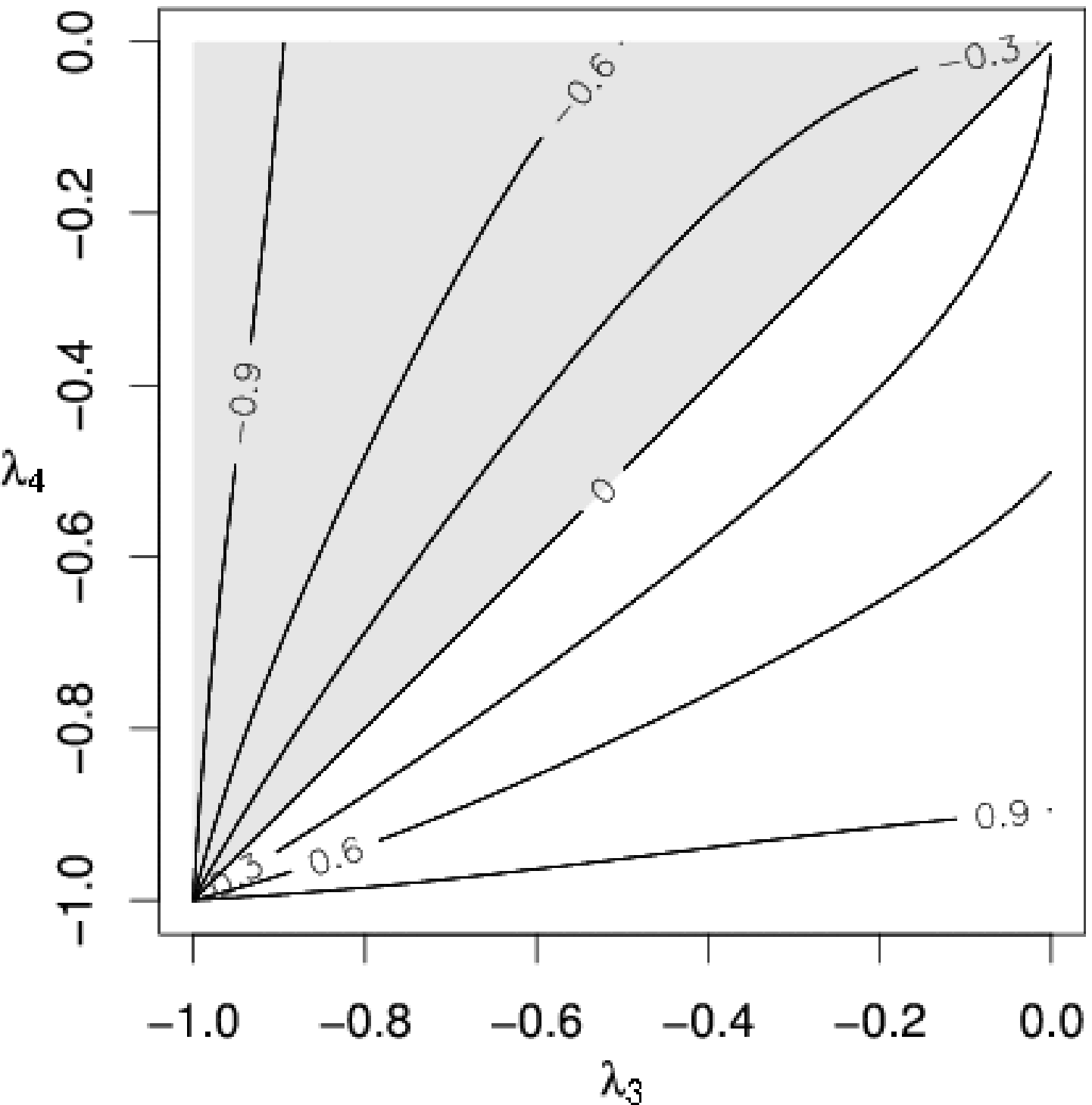}}
\subfigure[Contours of L-kurtosis $\tau_{4}$]{
 \includegraphics[width=0.48\columnwidth]{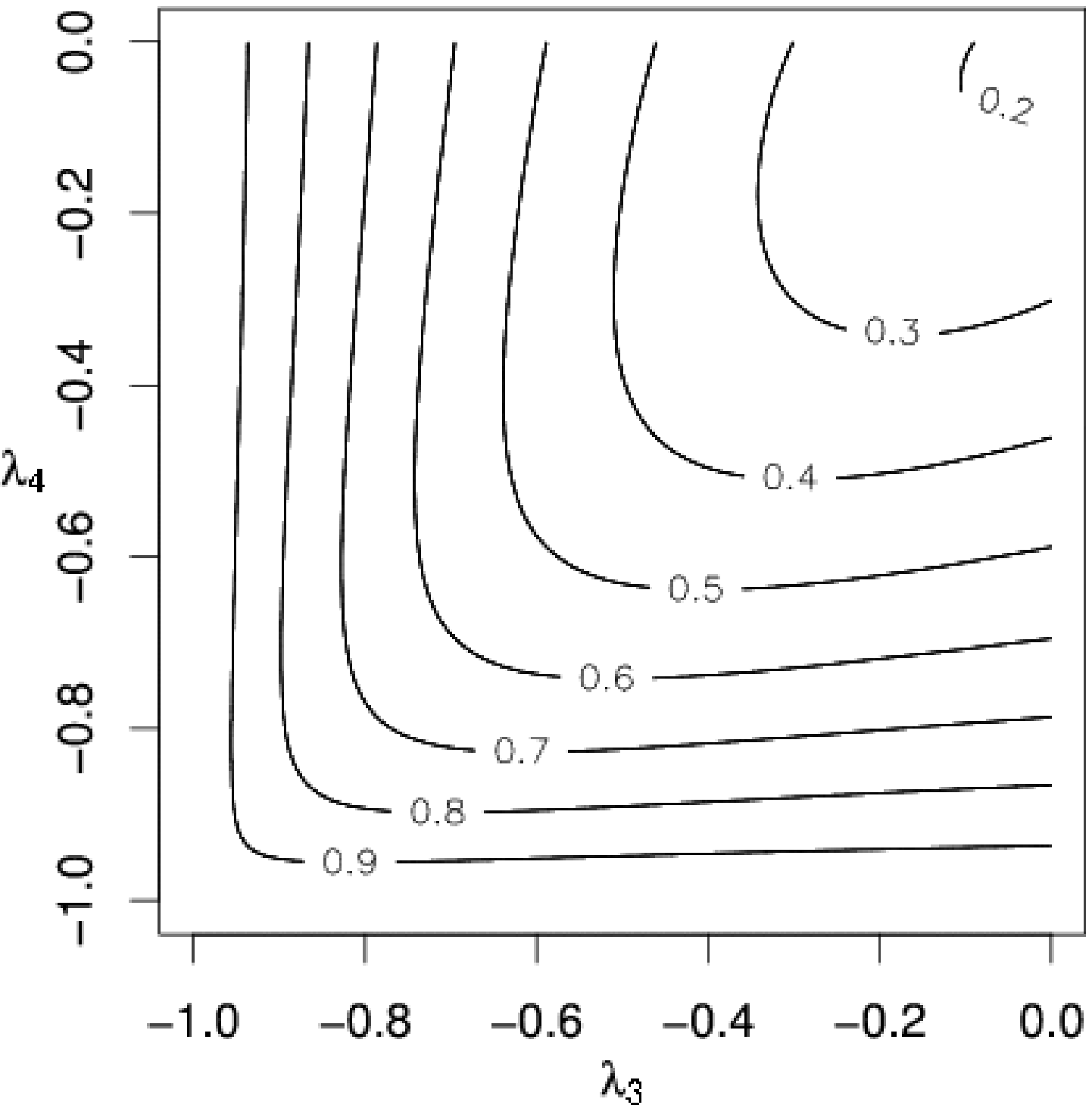}}
\end{center}
\caption{Contours of  $\tau_{3}$ and $\tau_{4}$ in region 4. The shading indicates negative values of $\tau_{3}$} \label{fig:reg4contours}
\end{figure}

\begin{figure}[htb]
\begin{center}
\subfigure[Contours of L-skewness $\tau_{3}$]{
 \includegraphics[width=0.48\columnwidth]{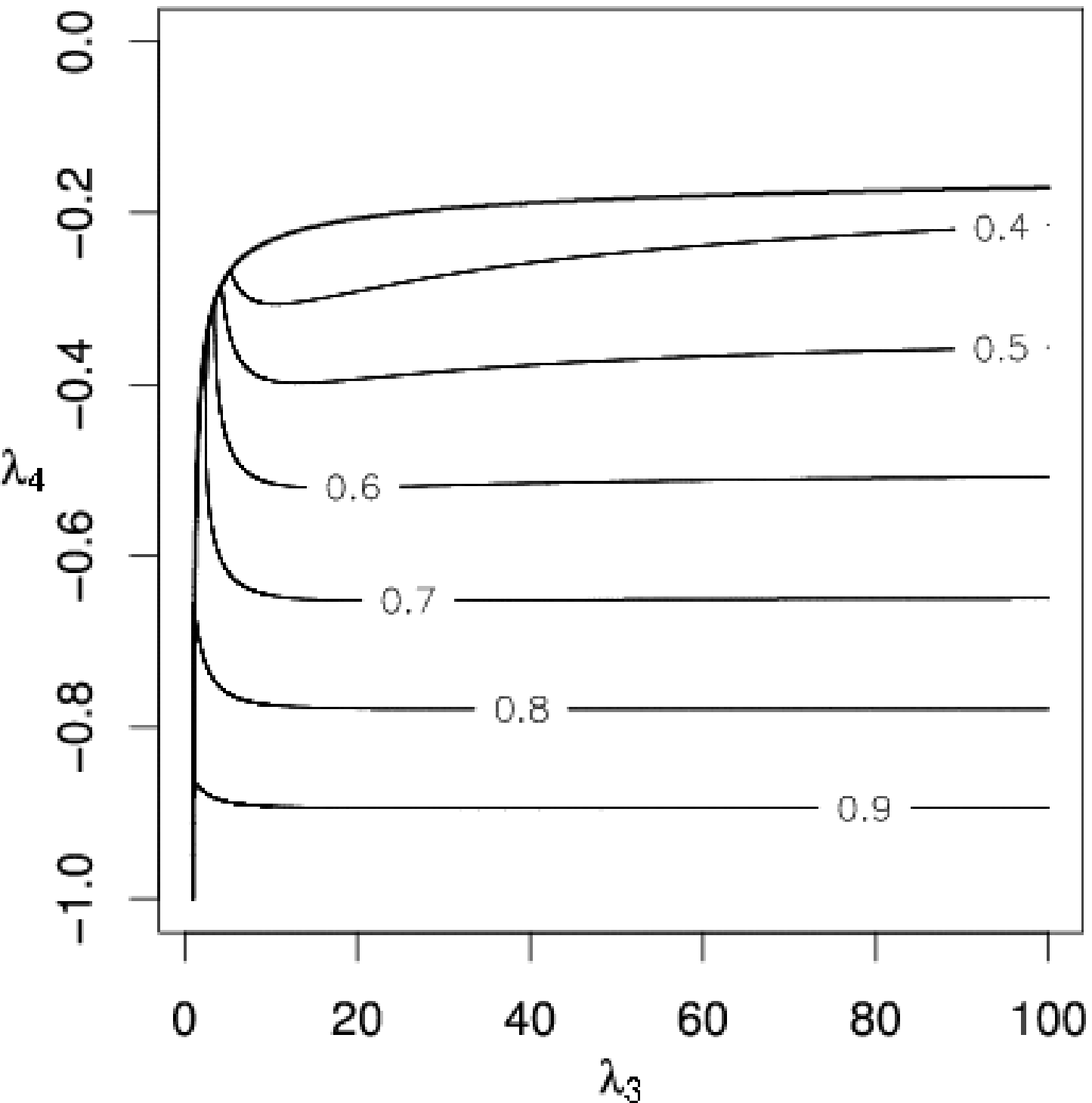}}
\subfigure[Contours of L-kurtosis $\tau_{4}$]{
 \includegraphics[width=0.48\columnwidth]{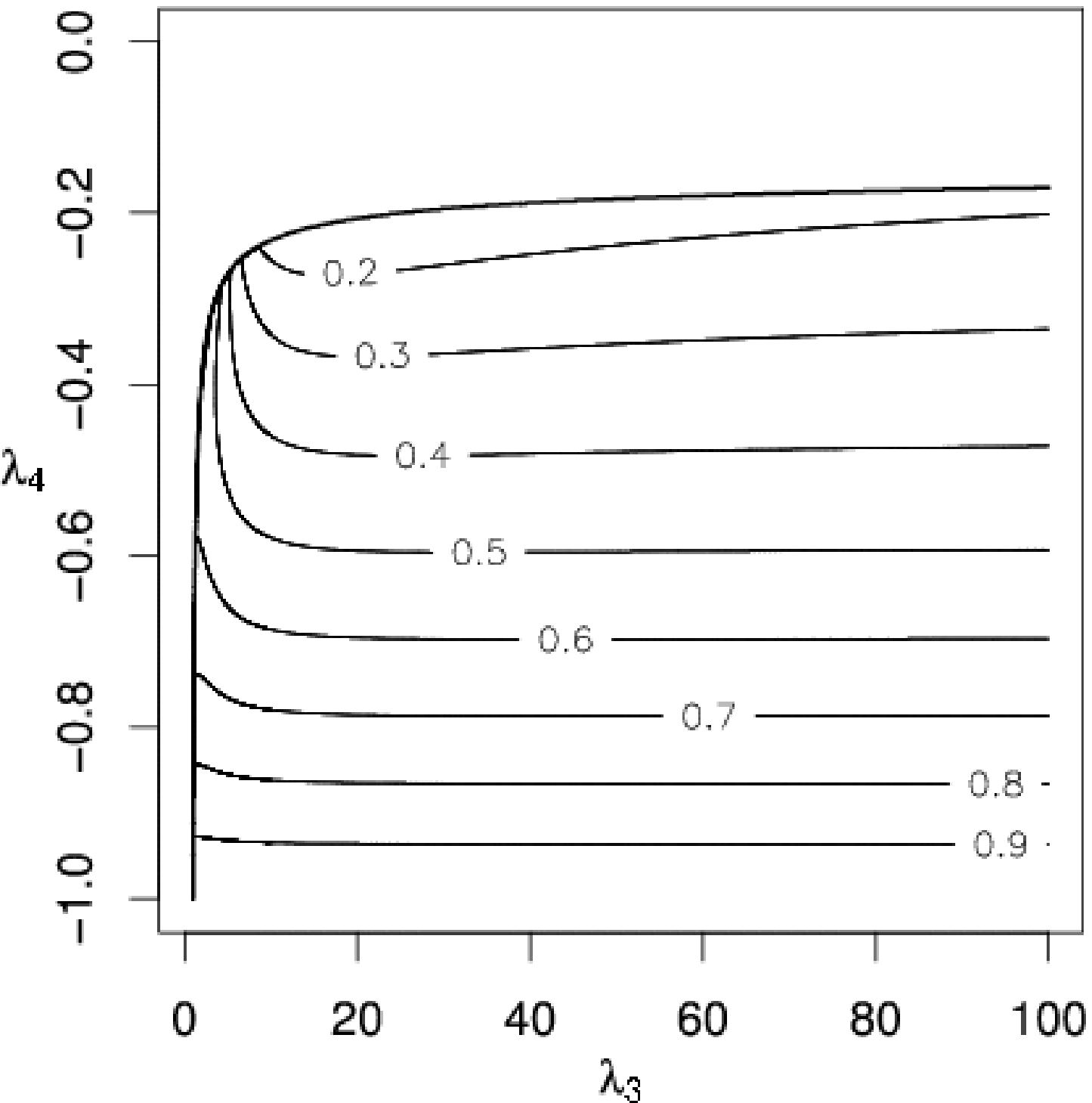}}
\end{center}
\caption{Contours of  $\tau_{3}$ and $\tau_{4}$ in region 6.} \label{fig:reg6contours}
\end{figure}

The contours of $\tau_{3}$ and $\tau_{4}$ in different GLD regions are shown in Figures~\ref{fig:reg3contours}, \ref{fig:reg4contours} and \ref{fig:reg6contours}. In region 3 the contours are rather complicated. From Figure~\ref{fig:reg3contours}(a) and ~\ref{fig:reg3contours}(c) it can be seen that there are three separate areas where $\tau_{3}$ has negative values. It is also seen that besides the line $\lambda_{3}=\lambda_{4}$, there are two curves where $\tau_{3}=0$. Naturally, the distributions defined by the curves are not symmetric; they just have zero L-skewness.

\section{Estimation by method of L-moments} \label{sec:estimation}
In method of L-moments, we first calculate L-moments $\hat{L}_{1}$, $\hat{L}_{2}$, $\hat{\tau}_{3}$ and $\hat{\tau}_{4}$ from the data. Then we numerically find parameters $\hat{\lambda}_{3}$ and $\hat{\lambda}_{4}$ that minimize some objective function that measures the distance between the L-skewness and L-kurtosis of the data $(\hat{\tau}_{3}, \hat{\tau}_{4})$ and the L-skewness and L-kurtosis of the estimated model $(\tau_{3}(\hat{\lambda}_{3},\hat{\lambda}_{4}),\tau_{4}(\hat{\lambda}_{3},\hat{\lambda}_{4}))$.
After that estimates $\hat{\lambda}_{1}$ and $\hat{\lambda}_{2}$ can be solved from equations~\eref{eq:L1} and \eref{eq:L2}. A natural objective function is the sum of squared distances
\begin{equation} \label{eq:sqcrit}
\big(\hat{\tau}_{3}-\tau_{3}(\hat{\lambda}_{3},\hat{\lambda}_{4})\big)^2+
\big(\hat{\tau}_{4}-\tau_{4}(\hat{\lambda}_{3},\hat{\lambda}_{4})\big)^2,
\end{equation}
which was used also by \citet{Asquith:lmomgld}.

\begin{figure}[htb]
\begin{center}
 \subfigure[Contours and the four solutions]{
 \includegraphics[width=0.48\columnwidth]{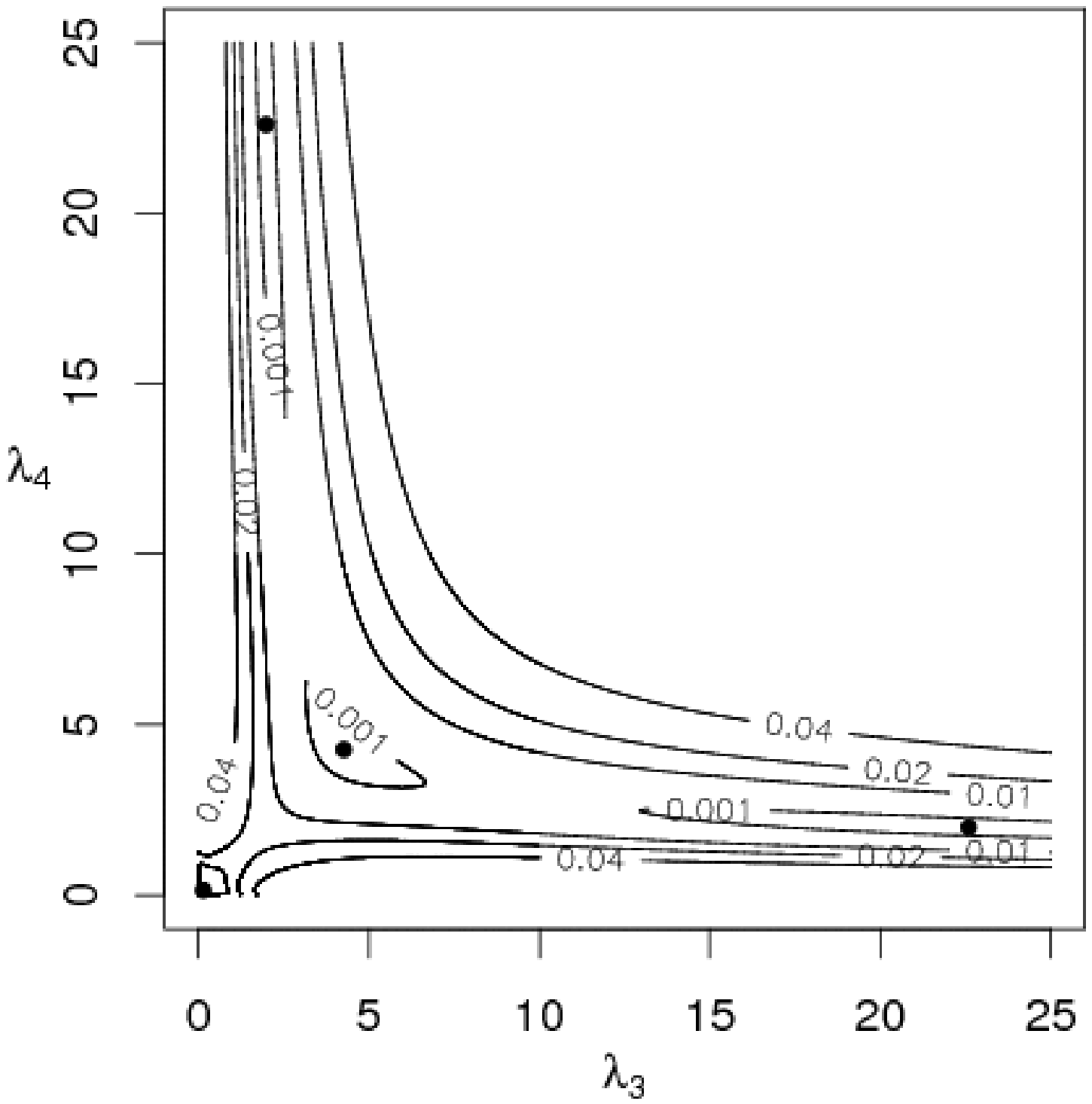}}
\subfigure[Contours in the neighborhood of $GLD(0,0.19,0.14,0.14)$]{
 \includegraphics[width=0.48\columnwidth]{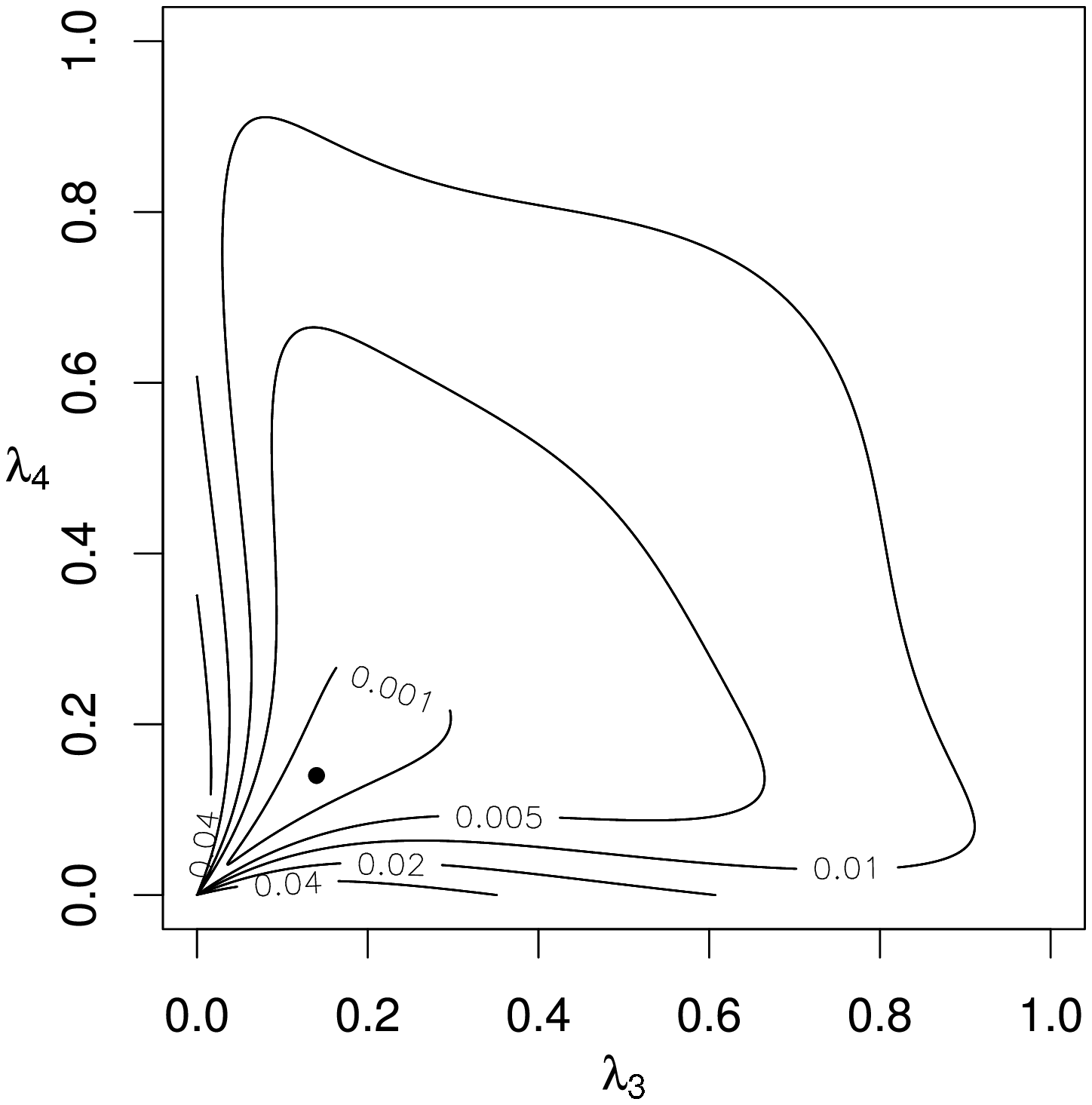}}
 \
\end{center}
\caption{Contours of objective function~\eref{eq:sqcrit} when $\hat{\tau}_{3}=0$ and $\hat{\tau}_{4}\approx 0.12305$. The dots show the two symmetric and the two asymmetric solutions.} \label{fig:crit_contours}
\end{figure}

We study the same simulation example that was studied by \citet{Fournier:estimating}. They compare five estimators in situation where a sample of 1000 observations is generated from $GLD(0,0.19,0.14,0.14)$, which is a symmetric distribution close to the standard normal distribution. We apply method of L-moments to the same problem and present the results in a comparable form. Before carrying out the simulations, we analyze the example using the results given in the preceding sections. Because $\lambda_{3}=\lambda_{4}=0.14$, $GLD(0,0.19,0.14,0.14)$ is a symmetric distribution from region 3. The theoretical L-moments $L_{1}=0$, $L_{2} \approx 0.60407$, $\tau_{3}=0$ and $\tau_{4} \approx 0.12305$ of $GLD(0,0.19,0.14,0.14)$ are obtained from equations~\eref{eq:L1}--\eref{eq:L4}. Because $\tau_{3}=0$ we may apply equation~\eref{eq:symmsolution} that gives the ``correct'' solution $\lambda_{3}=\lambda_{4}=0.14$ and another solution  $\lambda_{3}=\lambda_{4} \approx 4.26316$ also from region 3. By inspecting the contours in Figure~\ref{fig:reg3contours} we find that there exist also two asymmetric GLDs from region 3 with $\tau_{3}=0$ and $\tau_{4} \approx 0.12305$. The numerical minimization of objective function~\eref{eq:sqcrit} with suitable initial values gives solutions $(\lambda_{3} \approx 1.98, \lambda_{4} \approx 22.59)$ and $(\lambda_{3} \approx  22.59, \lambda_{4} \approx 1.98)$. Figure~\ref{fig:crit_contours} shows the contours of objective function~\eref{eq:sqcrit} when $\hat{\tau}_{3}=0$ and $\hat{\tau}_{4}\approx 0.12305$. It can be clearly seen that the achieved numerical solution depends on the initial values of $\lambda_{3}$ and $\lambda_{4}$.

In the simulation we choose a strategy where $\hat{\tau}_{4}$ computed from the data and the two solutions of equation~\eref{eq:symmsolution} are tried as initial values of $\lambda_{3}$ and $\lambda_{4}$. The numerical optimization of the objective function~\eref{eq:sqcrit} is carried out by the Nelder-Mead method \citep{NelderMead} in R \citep{R} and the computer used for the simulation is comparable to the computer used by \citet{Fournier:estimating}. The results from 10000 data sets of 1000 observations generated from $GLD(0,0.19,0.14,0.14)$ are presented in Table~\ref{tab:simuexample}. The first column of the table reports the L-moment estimates corresponding to the smaller initial solution which leads to the estimates corresponding to the generating GLD. The last three columns are copied from \citep{Fournier:estimating} and are based on 100 data sets of 1000 observations generated from the same GLD. The L-moment estimates of $\lambda_{1}$, $\lambda_{3}$ and $\lambda_{4}$ are unbiased and $\lambda_{2}$ is slightly underestimated. Method of L-moments is about 100 times faster than the Fournier et al. method and over 7000 times faster than the starship methods. In terms of the Kolmogorov-Smirnov statistics $E_{KS}$ between the observed data and the fitted model, method of L-moments gave slightly better fit than the starship methods with the used settings.

\begin{table}
\caption{Comparison of estimation methods. The mean and the estimated standard error of the L-moment estimators are computed from 1000 data sets generated from $GLD(0,0.19,0.14,0.14)$. The results for the other methods are based 100 generated data sets from the same GLD and are copied from \citep{Fournier:estimating}. Fournier et al. method is based on the percentile method and the minimization of the Kolmogorov-Smirnov distance in the $(\lambda_{3},\lambda_{4})$ space. Starship KS and starship AD methods correspond to the starship method \citep{King:starsip} using the Kolmogorov-Smirnov or the Anderson-Darling distance, respectively. \label{tab:simuexample}}
\begin{center}
\footnotesize
\begin{tabular}{llcccc}
Quantity & & L-moments & Fournier et al.  & Starship KS  & Starship AD\\
$\lambda_{1}$ & Mean &          0.00153  &    -0.00030  &  0.25010  &  0.00720 \\
              & Std error &0.10287  &   0.05010  &  0.10000  &  0.07300 \\
$\lambda_{2}$ & Mean &           0.18795  &   0.19960  &  0.20020  &  0.19400 \\
              & Std error &0.03625  &   0.05650  &  0.04760  &  0.03020 \\
$\lambda_{3}$ & Mean &           0.14012  &   0.15110  &  0.15530  &  0.14600 \\
              & Std error &0.03602  &   0.05150  &  0.04630  &  0.03170 \\
$\lambda_{4}$ & Mean &           0.13981  &   0.15060  &  0.14750  &  0.14350 \\
              & Std error &0.03573  &   0.05030  &  0.04510  &  0.02860 \\
Time (s)     & Mean &           0.02767  &   3.19010  &  206.99000  &  213.47000 \\
              & Std error &0.00575  &   0.15020  &  40.42010  &  42.15210 \\
$E_{KS}$     & Mean &           0.01597  &   0.01850  &  0.01680  &  0.01856 \\
              & Std error &0.00319  &   0.00430  &  0.00270  &  0.00360
\end{tabular}
\end{center}
\end{table}

\section{Conclusion} \label{sec:conclusion}
We have presented the L-moments of the GLD up to an arbitrary order and studied which values of L-skewness $\tau_{3}$ and L-kurtosis $\tau_{4}$ can be achieved by the GLD. For the symmetric case the boundaries were derived analytically and in the general case numerical methods were used. It was found that with an exception of the smallest values of $\tau_{4}$, the GLD covers all possible $(\tau_{3},\tau_{4})$ pairs and often there are two or more distributions sharing the same $\tau_{3}$ and $\tau_{4}$. The example in Section~\ref{sec:boundaries} demonstrates a situation where there are four GLD members sharing the same $\tau_{3}$ and $\tau_{4}$.

We argue that L-moments are natural descriptive statistics for the GLD because they, differently from the central moments, can be expressed in closed form. The existence of the closed form presentation of the GLD L-moments follows from the fact that the inverse distribution function of the GLD is available in closed form and can be integrated. The relation between L-moments and distributions defined by the inverse distribution function is not restricted only to the GLD \citep{qm}. The existence of the first four central moments of the GLD requires that $\lambda_{3},\lambda_{4} > -1/4$ whereas the existence of L-moments of any order requires only that $\lambda_{3},\lambda_{4} > -1$. Thus, using the L-moments we can characterize wider subset of the GLD than using the central moments. \citet{Asquith:lmomgld} concluded that even wider subset of the GLD could be characterized using the trimmed L-moments \citep{Elamir:trimmedlmoments}.

The results presented in this paper can be utilized in model selection and estimation. The characterization by L-moments gives an insight into which $\tau_{3}$ and $\tau_{4}$ are available in each GLD region. For instance, there are symmetric GLD members available from both region 3 and 4, if $\tau_{4}>1/6$. This kind of information is useful when making the decision whether the GLD is a potential model for certain data. The results are also useful in the estimation of the GLD parameters. The parameters can be estimated directly by method of L-methods or the L-moment estimates can be used as starting values for other estimation methods. In both cases, we can use the L-moment ratio boundaries to specify the potential GLD regions where we should search for the parameter estimates. The choice between alternative solutions can be based on the type of the domain (bounded/unbounded) or some other additional criterion such as the Kolmogorov-Smirnov statistic or L-moment ratios $\tau_{5}$ and $\tau_{6}$.

In the simulation example, method of L-moments compared favorably to more complicated estimation methods. Estimation by method of L-moments was the fastest in the comparison and the estimates were unbiased or nearly unbiased. The goodness of fit measured by the Kolmogorov-Smirnov statistic was the same or better than with the alternative estimation methods. More simulations are needed to find out how general these results are.

\section*{Acknowledgement}
The authors thank William H. Asquith for useful comments.


\begin{thebibliography}{28}
\expandafter\ifx\csname natexlab\endcsname\relax\def\natexlab#1{#1}\fi
\expandafter\ifx\csname url\endcsname\relax
  \def\url#1{\texttt{#1}}\fi
\expandafter\ifx\csname urlprefix\endcsname\relax\def\urlprefix{URL }\fi

\bibitem[{Asquith(2007)}]{Asquith:lmomgld}
Asquith, W.~H., 2007. {L}-moments and {TL}-moments of the generalized lambda
  distribution. Computational Statistics \& Data Analysis 51, 4484--4496.

\bibitem[{Bergevin(1993)}]{Bergevin:analysis}
Bergevin, R.~J., 1993. An analysis of the generalized lambda distribution.
  Master's thesis, Air Force Institute of Technology, available from
  \url{http://stinet.dtic.mil}.

\bibitem[{Bigerelle et~al.(2005)Bigerelle, Najjar, Fournier, Rupin, and
  Iost}]{Bigerelle:fatigue}
Bigerelle, M., Najjar, D., Fournier, B., Rupin, N., Iost, A., 2005. Application
  of lambda distributions and bootstrap analysis to the prediction of fatigue
  lifetime and confidence intervals. International Journal of Fatigue 28,
  223--236.

\bibitem[{Corrado(2001)}]{Corrado:optionpricing}
Corrado, C.~J., 2001. Option pricing based on the generalized lambda
  distribution. Journal of Future Markets 21, 213--236.

\bibitem[{Elamir and Seheult(2003)}]{Elamir:trimmedlmoments}
Elamir, E.~A., Seheult, A.~H., 2003. Trimmed {L}-moments. Computational
  Statistics \& Data Analysis 43, 299--314.

\bibitem[{Fournier et~al.(2007)Fournier, Rupin, Bigerellle, Najjar, Iost, and
  Wilcox}]{Fournier:estimating}
Fournier, B., Rupin, N., Bigerellle, M., Najjar, D., Iost, A., Wilcox, R.,
  2007. Estimating the parameters of a generalized lambda distributions.
  Computational Statistics \& Data Analysis 51, 2813--2835.

\bibitem[{Headrick and Mugdadib(2006)}]{Headrick:simulating}
Headrick, T.~C., Mugdadib, A., 2006. On simulating multivariate non-normal
  distributions from the generalized lambda distribution. Computational
  Statistics \& Data Analysis 50, 3343--3353.

\bibitem[{Hosking(1990)}]{Hosking:Lmoments}
Hosking, J., 1990. L-moments: Analysis and estimation of distributions using
  linear combinations of order statistics. Journal of Royal Statistical Society
  B 52~(1), 105--124.

\bibitem[{Hosking(2006)}]{Hosking:characterization}
Hosking, J. R.~M., 2006. On the characterization of distributions by their
  {L}-moments. Journal of Statistical Planning and Inference 136~(1), 193--198.

\bibitem[{Jones(2004)}]{Jones:someexpressions}
Jones, M.~C., 2004. On some expressions for variance, covariance, skewness and
  {L}-moments. Journal of Statistical Planning and Inference 126, 97--106.

\bibitem[{Karian and Dudewicz(1999)}]{egld:percentile}
Karian, Z.~A., Dudewicz, E.~J., 1999. Fitting the generalized lambda
  distribution to data: A method based on percentiles. Communications in
  Statistics: Simulation and Computation 28~(3), 793--819.

\bibitem[{Karian and Dudewicz(2000)}]{KarianDudewicz:gldbook}
Karian, Z.~A., Dudewicz, E.~J., 2000. Fitting Statistical Distributions: The
  Generalized Lambda Distribution and Generalized Bootstrap Methods. Chapman \&
  Hall/CRC Press, Boca Raton, Florida.

\bibitem[{Karian and Dudewicz(2003)}]{Karian:comparison}
Karian, Z.~A., Dudewicz, E.~J., 2003. Comparison of {GLD} fitting methods:
  {S}uperiority of percentile fits to moments in {$L^{2}$} norm. Journal of
  Iranian Statistical Society 2~(2), 171--187.

\bibitem[{Karian et~al.(1996)Karian, Dudewicz, and McDonald}]{egld:theory}
Karian, Z.~A., Dudewicz, E.~J., McDonald, P., 1996. The extended generalized
  lambda distribution system for fitting distributions to data: History,
  completion of theory, tables, applications, the ``final word'' on moment
  fits. Communications in Statistics: Simulation and Computation 25~(3),
  611--642.

\bibitem[{Karvanen(2003)}]{gldrnd}
Karvanen, J., 2003. Generation of correlated non-{G}aussian random variables
  from independent components. In: Proceedings of Fourth International
  Symposium on Independent Component Analysis and Blind Signal Separation,
  ICA2003. pp. 769--774.

\bibitem[{Karvanen(2006)}]{qm}
Karvanen, J., 2006. Estimation of quantile mixtures via {L}-moments and trimmed
  {L}-moments. Computational Statistics \& Data Analysis 51~(2), 947--959.

\bibitem[{Karvanen et~al.(2002)Karvanen, Eriksson, and Koivunen}]{vlsi}
Karvanen, J., Eriksson, J., Koivunen, V., 2002. Adaptive score functions for
  maximum likelihood {ICA}. Journal of VLSI Signal Processing 32, 83--92.

\bibitem[{King(2006)}]{R_gld}
King, R., 2006. gld: Basic functions for the generalised (Tukey) lambda
  distribution. R package version 1.8.1.

\bibitem[{King and MacGillivray(1999)}]{King:starsip}
King, R., MacGillivray, H., 1999. A startship estimation method for the
  generalized lambda distributions. Australian and New Zealand Journal of
  Statistics 41~(3), 353--374.

\bibitem[{Lakhany and Mausser(2000)}]{Lakhany:estimating}
Lakhany, A., Mausser, H., 2000. Estimating parameters of generalized lambda
  distribution. Algo Research Quarterly 3~(3), 47--58.

\bibitem[{Lampasi et~al.(2005)Lampasi, {Di Nicola}, and
  Podest{\`a}}]{Lampasi:glduncertainty}
Lampasi, D.~A., {Di Nicola}, F., Podest{\`a}, L., 2005. The generalized lambda
  distribution for the expression of measurement uncertainty. In: Proceedings
  of IMTC--2005 IEEE Instrumentation and Measurement Technology Conference. pp.
  2118--2133.

\bibitem[{Nelder and Mead(1965)}]{NelderMead}
Nelder, J.~A., Mead, R., 1965. A simplex algorithm for function minimization.
  Computer Journal 7, 308--–313.

\bibitem[{{\"O}zt{\"u}rk and Dale(1985)}]{Ozturk:leastsquares}
{\"O}zt{\"u}rk, A., Dale, R., 1985. Least squares estimation of the parameters
  of the generalized lambda distribution. Technometrics 27~(1), 81--84.

\bibitem[{Pal(2005)}]{Pal:processgld}
Pal, S., 2005. Evaluation of non-normal process capability indices using
  generalized lambda distributions. Quality Engineering, 77--85.

\bibitem[{{R Development Core Team}(2006)}]{R}
{R Development Core Team}, 2006. R: A Language and Environment for Statistical
  Computing. R Foundation for Statistical Computing, Vienna, Austria.
\newline\urlprefix\url{http://www.R-project.org}

\bibitem[{Ramberg and Schmeiser(1974)}]{Ramberg:74}
Ramberg, J.~S., Schmeiser, B.~W., 1974. An approximate method for generating
  asymmetric random variables. Communications of the ACM 17, 78--82.

\bibitem[{Schnute et~al.(2006)Schnute, Boers, Haigh, and {and
  others}}]{R_PBSmapping}
Schnute, J., Boers, N., Haigh, R., {and others}, 2006. PBSmapping: PBS Mapping
  2. R package version 2.09.

\bibitem[{Su(2007)}]{Su:numericalmaximum}
Su, S., 2007. Numerical maximum log likelihood estimation for generalized
  lambda distributions. Computational Statistics \& Data Analysis 51,
  3983--3998.

\end{thebibliography}
\end{document}